%% file: 2000-3.tex
\documentclass{gtart} 

\usepackage{rlepsf, amsmath}
\let\relabela\adjustrelabel

\input gtoutput
\volumenumber{4}\papernumber{3}\volumeyear{2000}
\pagenumbers{117}{148}\published{29 February 2000}
\proposed{David Gabai}\seconded{Jean-Pierre Otal, Walter Neumann}
\received{16 July 1999}\revised{9 November 1999}
\accepted{20 February 2000}

\input{amsmacs}  

\long\def\state#1#2{
\medskip\par\noindent
{\bf #1}\qua 
{\sl #2}
\par\medskip
}

\renewcommand{\bfheading}[1]{\par\medskip\noindent {\bf #1}\qua\ignorespaces}

\let\segment\sh
\def\S{Section }

\newcommand{\transL}{L_3}
\newcommand{\projmap}{{\mathbf p}}
\newcommand{\proj}{{\mathbf P}}
\newcommand{\inj}{\operatorname{inj}}

\begin{document}

\title{Kleinian groups and the complex of curves}
\author{Yair N Minsky}

\address{Department of Mathematics, SUNY at Stony Brook\\Stony Brook, 
NY 11794, USA}

\email{yair@math.sunysb.edu}

\begin{abstract}
We examine the internal geometry of a Kleinian surface group and its
relations to the asymptotic geometry of its ends, using the
combinatorial structure of the complex of curves on the surface.
Our main results give necessary conditions for the Kleinian group to have
`bounded geometry' (lower bounds on injectivity radius) in terms of a
sequence of coefficients (subsurface projections) computed using the
ending invariants of the group and the complex of curves.

These results are directly analogous to those obtained in the case of
punctured-torus surface groups. In that setting the ending invariants
are points in the closed unit disk and the coefficients are closely
related to classical continued-fraction coefficients. The estimates
obtained play an essential role in the solution of Thurston's ending
lamination conjecture in that case.
\end{abstract}




\maketitlepage

\section{Introduction}\label{intro}

In this paper we examine some new connections between the internal geometry of
hyperbolic 3--manifolds and the asymptotic geometry of their ends. 
Our main theorem is a necessary condition for ``bounded geometry''
(lower bounds on injectivity radius), in terms of combinatorial
properties of the asymptotic geometry, which we think of as analogous
to bounds on continued fraction expansions. We conjecture also that
this necessary condition is sufficient.

The motivating problem for this analysis is Thurston's Ending
Lamination Conjecture, which states that a hyperbolic 3--manifold with
finitely-generated fundamental group is 
determined, up to isometry, by its topological type and a collection of
invariants that specify the asymptotic geometry of its ends.
Even in the geometrically finite case, where Thurston's
conjecture reduces to classical work of Ahlfors--Bers et al,
the theory is not geometrically explicit:
direct estimates of internal geometry from the end invariants are quite
hard to come by.

We will approach this problem by studying 
the combinatorial/geometric structure of 
an object known as  a {\em complex of curves on a surface},  a simplicial
complex $\CC(S)$ whose vertices are homotopy classes of simple closed
curves on a surface $S$ (see \S\ref{complex} for complete definitions).
Adjacency relations in these complexes are related to elementary
homotopies between special types of surfaces in hyperbolic
3--manifolds, in a way that gives some control of the internal
geometry. Our main methods of proof in this paper depend on some
important properties of pleated surfaces in hyperbolic 3--manifolds,
originally discovered and applied by Thurston.
The results we obtain here should be viewed as the first
steps of a project to exploit these connections.

\segment{Statements of results}

In this paper, 
a {\em Kleinian surface group} is a discrete, faithful, type-preserving
representation $\rho\co \pi_1(S) \to \PSL 2(\C)$, where $S = S_{g,p}$
is an oriented 
genus $g$ surface with $p$ punctures, and the
type-preserving condition means that the image of a loop around any
puncture is parabolic. We exclude $S_{0,p}$ for $p\le 2$ and
$S_{1,0}$, for which there are only elementary representations. 
Let $N_\rho$ be the quotient $\Hyp^3/\rho(\pi_1(S))$.
Such representations arise naturally in the
theory when considering hyperbolic 3--manifolds homotopy-equivalent to
manifolds with incompressible boundary, and restricting to the
boundary groups.

If the group is {\em quasi-Fuchsian} then
its action on the Riemann sphere has two invariant disks whose
quotients give two Riemann surface structures on $S$, or two points in
the Teichm\"uller space $\TT(S)$, which we label
$\nu_-(\rho)$ and $\nu_+(\rho)$. The work of Ahlfors--Bers
  \cite{ahlfors-bers,bers:simultaneous} shows that this gives a
  parametrization of all quasi-Fuchsian representations up to
  isometry by
  $\TT(S)\times\TT(S)$, which is in fact a holomorphic
  isomorphism. And yet, it is hard to answer simple questions such as, 
given $\nu_\pm$, what is the length of the shortest geodesic (or any
given geodesic) in the manifold $N_\rho$?

For general representations, Thurston and Bonahon
\cite{wpt:notes,bonahon} showed how to generalize $\nu_\pm$ using {\em
  ending laminations} (see \S\ref{proof of A}). The same questions are
difficult to answer in this case, and in addition
Thurston's Ending Lamination Conjecture (the analogue of the
Ahlfors--Bers parameterization), which in this setting states that 
$\nu_\pm$ determine $\rho$ up to isometry, is open in most cases.

Our first theorem gives a
necessary condition, in terms of $\nu_\pm(\rho)$, 
for $\rho$ to have {\em  bounded geometry } in the following sense.
Define the {\em non-cuspidal injectivity radius} $\inj_0(\rho)$ to be
half the infimum of the lengths of all closed geodesics in $N_\rho$.
We say $\rho$ has bounded geometry if
$\inj_0(\rho)>0$. 

For any essential subsurface $Y\subset S$
we will define a ``projection'' $\pi_Y$ that takes the invariants
$\nu_\pm$ to 
elements of the complex of 
curves of $Y$, and we denote
the distance measured in this complex
by $d_Y(\nu_-,\nu_+)$.

\state{Theorem A}{
For any Kleinian surface group $\rho$ with ending invariants
$\nu_\pm$,  if
$$
\sup_{Y\subset S} d_Y(\nu_+,\nu_-) = \infty
$$
then
$$
\inj_0(\rho) = 0,
$$
where the supremum is over proper essential subsurfaces $Y$ in $S$
not all of whose boundaries are mapped to parabolics.
}

See below for some examples of surface groups to which this theorem
applies, and some examples supporting the converse implication.

Theorem A is a consequence of Theorem B, in which we
consider the set of geodesics in $N_\rho$ satisfying a
certain length bound. If $\alpha$ is a finite set of homotopy classes
of simple closed
curves in $S$, let $\ell_\rho(\alpha)$ denote the total length of
their geodesic representatives in $N_\rho$.
For $L>0$, let
$\CC_0(\rho,L)$ denote the set of homotopy classes $\alpha$ of simple closed
curves in $S$ with $\ell_\rho(\alpha) \le L$. An understanding of
this set is
crucial to characterizing the internal geometry of $N_\rho$.
The projection $\pi_Y$ is defined on this set and the diameter of the
image is denoted 
$\diam_Y(\CC_0(\rho,L))$. We will prove:

\bfheading{Theorem B}{\sl
Given a surface $S$, $\ep>0$ and $L>0$, there exists $K$ 
so that, if $\rho\co \pi_1(S)\to \PSL 2(\C)$ is a Kleinian surface group
and $Y$ is a proper essential subsurface of $S$, then 
$$
\diam_Y(\CC_0(\rho,L)) \ge K
\quad\implies\quad
\ell_\rho(\boundary Y) \le \ep.
$$
}

The relation between the invariants $\nu_\pm(\rho)$ and the set of bounded
curves $\CC_0(\rho,L)$ is, for example in the case where
$\nu_\pm$ are laminations rather than points in the Teichm\"uller
space, that  the elements of $\CC_0(\rho,L)$ for suitable $L$ will
accumulate exactly on $\nu_\pm$ in the space of laminations on
$S$. This gives the connection between Theorem B and Theorem A.

\segment{Examples}

Let us sketch some examples of surface groups where the hypothesis of Theorem
A holds. The basic idea is similar to that in Bonahon--Otal \cite{bonahon-otal}
and Thurston \cite{wpt:II}, where surface groups $\rho\co \pi_1(S)\to
\PSL 2(\C)$ are constructed with arbitrarily short geodesics. In those
examples the projections $\pi_Y(\nu_-,\nu_+)$ which are large
correspond to annuli $Y$ whose cores are the short curves.
We leave out many details, since they are not needed for the proofs of
the main theorems, but the reader is referred to 
\cite{wpt:surfaces,travaux,casson}
for details about laminations and pseudo-Anosov maps, and 
\cite{mcmullen:renormbook}
for details about quasi-Fuchsian groups and Bers slices.

Let $Y_1,Y_2,\ldots$ be a sequence of essential subsurfaces of $S$ and
let $\Psi_i\co S\to S$ be homeomorphisms so that 
$\Psi_i = id$ outside $Y_i$, and $\Psi_i|_{Y_i}$ are pseudo-Anosov.
For simplicity it is probably a good idea to take the $Y_i$ and
$\Psi_i$ from a finite set, eg, $Y_i = Y_{i+2}$ 
and  $\Psi_i = \Psi_{i+2}$ for each $i$. 
Assume also that $\boundary Y_i$ and $\boundary Y_{i+1}$ intersect
each other essentially.

Let $\{m_i\}$ be an increasing sequence of positive integers.
Define
$$
h_j^k = \Psi_j^{m_j} \circ \cdots \circ \Psi_{k-1}^{m_{k-1}}
$$
for integers $0<j<k$, and extend to all pairs of integers $j,k>0$ by requiring
$h_k^j = (h_j^k)^{-1}$. Note that $h_j^k\circ h_k^m = h_j^m$.
Now fix a conformal structure $\sigma$ on $S$ and consider the
sequence of quasi-Fuchsian groups
$\rho_k = QF(\sigma,h_1^k(\sigma))$ as $k\to \infty$.
(The notation $QF(\sigma,\tau)$ refers to a quasi-Fuchsian group
whose Riemann surface structures at infinity are $\nu_-=\sigma$,
$\nu_+=\tau$. This is determined up to conjugacy in $\PSL 2(\C)$.)

Let $Z_j = h_1^j(Y_j)$. We claim that, with an {\it a priori} choice of
sufficiently large numbers $m_j$, the coefficients
$c_j(k) = d_{Z_j}(\sigma, h_1^k(\sigma))$ can be made as large as we
please (where $j<k$).
We refer the reader to sections \ref{complex} and \ref{proof of A}
where these coefficients are defined, but note that in 
the case of $d_Z(\sigma,\tau)$ where $\sigma$ and $\tau$ are conformal
structures, and $Z$ is not an annulus, we simply
consider minimal-length curves $\alpha$ in $\sigma$ and $\beta$ in
$\tau$
which intersect $Z$, and measure the
``elementary
move'' distance in the complex $\CC'(Z)$  between the arc systems 
$\alpha\intersect Z$ and $\beta\intersect Z$.

Applying $h_j^1$ to $Z_j$, $\sigma$ and $h_1^k(\sigma)$, we 
obtain
$ c_j(k) = d_{Y_j}(h_j^1(\sigma), h_j^k(\sigma))$.
Let $\lambda_{i}^+$ and $\lambda_i^-$ be the stable and unstable
laminations, respectively, of $\Psi_i$.
Then   choosing $m_{j-1}$ sufficiently large we can guarantee that
$h_j^1(\sigma)$ is as close as we like (in Thurston's compactification
of Teichm\"uller space, say) to $\lambda_{j-1}^-$, and choosing $m_j$
sufficiently large $h_j^k(\sigma)$ is as close as we like
to $\lambda_{j}^+$. Note that our assumption that $Y_j$ and $\Psi_j$
vary in a finite set makes it easier to verify this in a uniform way.
Now since $\lambda_{j-1}^-$ intersects $Y_j$ in some finite number of
arcs up to homotopy, whereas $\lambda_j^+$ is a lamination that fills
$Y_j$, it is not hard to see (for example see
\cite{masur-minsky:complex1}) that the combinatorial distance 
$d_{Y_j}(h_j^1(\sigma),h_j^k(\sigma))$ can be made as large as we like
by increasing $m_j$ --- independently of $k$ as long as $k>j$.

Let $c_j = \lim_{k\to\infty} c_j(k)$. We can restrict to a subsequence
so that these limits exist, and by choosing $m_j\to\infty$ as above,
we guarantee that $c_j\to\infty$.
Now
considering the limit (in the Bers slice) of $\rho_k$, we obtain a
limiting representation $\rho$ with invariants $\nu_-=\sigma$, $\nu_+
= \lim_{k\to \infty} h_1^k(\sigma)$, and it is not hard to see that
$d_{Z_j}(\nu_-,\nu_+)$ are estimated (to within a bounded distance)
by $c_j$.
Theorem A implies that there are arbitrarily short geodesics in this
manifold. In fact Theorem B tells us that the boundary curves
$\boundary Z_j$ have lengths going to 0 as $j\to\infty$. The estimates
of Theorem B hold for all the groups $\rho_k$ in the sequence.


In the opposite direction (supporting the conjectural converse of
Theorem A), consider the representations
$\rho\co \pi_1(S)\to \PSL 2(\C)$ that arise as fibres of hyperbolic
manifolds that fibre over the circle (Thurston \cite{wpt:II}).
In this case the end-invariants are the stable and unstable
laminations $\lambda^\pm$ of the pseudo-Anosov monodromy $\Psi$ of the
fibration. Such a manifold certainly has bounded geometry since it 
covers a compact manifold, and we can show
that the distances $d_Y(\nu_-,\nu_+)$
are bounded over all subsurfaces $Y$: Let $\PML(S)$ denote the
projective measured lamination space \cite{penner-harer}, and 
let $K$ be any compact subset of $\PML(S) \setminus \{\nu_-,\nu_+\}$.
Fix a hyperbolic metric on $S$, and realize all laminations
geodesically in this metric.
For any $\lambda\in K$, $(S-\lambda)\intersect \nu_+$ is a collection
of arcs with length bounded by some finite $L(\lambda)$ (since $\nu_+$
fills $S$). Furthermore $L$ is bounded on
$K$: if not, then by compactness of $K$ there is some convergent
sequence $\lambda_i\to \lambda\in K$ for which
$L(\lambda_i)\to\infty$, and hence there is a sequence of arcs
$\alpha_i$ of $\nu_+$ of length going to infinity with
$\lambda_i$ disjoint from $\alpha_i$. Since $\nu_+$ is minimal these
arcs converge to $\nu_+$ and we conclude that $\nu_+$ and $\lambda$
have no transverse intersection, and hence are the same in $\PML(S)$
since $\nu_+$ 
fills and is uniquely ergodic. This contradicts the assumption that
$K$ misses $\nu_+$. The 
same argument can be made for $\nu_-$. Thus, for any subsurface $Y$
whose boundary $\boundary Y$, viewed as a lamination, determines an
element of $K$, there is a uniform upper bound $L_K$ on the lengths of 
components of $\nu_+\intersect Y$ and $\nu_- \intersect Y$. This 
gives, as in Section \ref{proof of B}, an upper bound $D_K$ on
$d_Y(\nu_-,\nu_+)$ over all $Y$ with $\boundary Y$ in $K$.
Now the action of the pseudo-Anosov $\Psi$ on
$\PML(S)-\{\nu_-,\nu_+\}$ has some compact fundamental domain $K$, 
and since $\nu_\pm$ are invariant under $\Psi$, the bound
$d_Y(\nu_-,\nu_+)\le D_K$ holds 
for $\boundary Y$ in each translate $\Psi^n(K)$. Thus the bound holds
for all $Y$.

A more complete class of examples is obtained when $S$ is a
once-punctured torus.
In this case, the set
$\CC_0(S)$ of simple closed curves in $S$ can be identified with
$\hat \Q = \Q\union\{\infty\}$ 
by taking each curve to
its slope in the homology $H_1(S) = \Z^2$ after fixing a basis.
The natural completion of $\hat\Q$ is the circle $\hat\R =
\R\union\{\infty\}$. 
The Teichm\"uller space $\TT(S)$ can 
be identified with the hyperbolic plane $\Hyp^2$, and the natural
compactification of $\Hyp^2$ by $\hat\R$ is exactly Thurston's
compactification of $\TT(S)$ by the  space of
projective measured laminations on $S$.
Thus $\nu_+$ and $\nu_-$ are points in the closed 2--disk $\Hyp^2\union\hat\R$.
The essential subsurfaces $Y$ in $S$ are all annuli, so they can 
again be parametrized by $\hat\Q$, by associating to each $Y$ the
slope of its core curve.
For simplicity, consider the special case where $\nu_- =
\infty$ and $\nu_+\in\R\setminus\Q$. 
Then in fact $d_Y(\nu_-,\nu_+)$ is uniformly bounded except when $Y$'s
slope is a continued-fraction approximant of $\nu_+$, and in that case
it is equal (up to bounded error) to the corresponding coefficient
$w(Y)$ in the continued-fraction expansion  of $\nu_+$.

Theorem B implies that $\ell_\rho(\boundary Y)$ is small if $w(Y)$
is large, and in fact there is an explicit estimate  for the {\em complex}
translation length $\lambda_\rho = \ell_\rho + i\theta_\rho$ (where
$\theta_\rho$ gives the rotation of the associated isometry), of the
form $\lambda_\rho(\boundary Y)\asymp 
(1/w(Y)^2 + 2\pi i/w(Y))$. See \cite{minsky:torus} for more details.

\segment{Outline of the paper}

Section \ref{complex} outlines the basic definitions and properties of
the complexes of curves $\CC(S)$. 
Section \ref{pleated surfaces} discusses some
properties of pleated surfaces in hyperbolic 3--manifolds,
concentrating on Thurston's Uniform Injectivity Theorem and some
consequences, old and new. A number of the proofs in these two sections are
sketched because they are reasonably well known, although
perhaps hard to find in print.

In Section \ref{proof of B} we give the proof of Theorem B. This is
the main construction of the paper, and the main tools for the proof
are our extension, Theorem \ref{Efficiency of pleated surfaces}, of
Thurston's  theorem on Efficiency of pleated surfaces, and Lemma
\ref{Short bridge arcs} on comparing pleated surfaces that share part of their
pleating locus.
The proof of Theorem A appears
next in Section \ref{proof of A}, after a short discussion of the end
invariants $\nu_\pm$.

\section{Complexes of curves and subsurface projections}
\label{complex}
We recall here the definitions of Harvey's complexes of curves
\cite{harvey:boundary} and related complexes from
\cite{masur-minsky:complex2}. 

Let $S=S_{g,p}$ be an orientable surface of genus $g$ with $p$
punctures.  We will always assume that $S$ is not a sphere with 0,1 or
3 holes, nor a torus with 0 holes. 
The case of the once-punctured torus $S_{1,1}$ and quadruply-punctured
sphere, $S_{0,4}$,  are special and will be
treated below, as will the case of the annulus $S_{0,2}$, which will 
only occur as a subsurface of larger surfaces.
We call all other cases ``generic''. 

If $S$ is generic, we define $\CC(S)$
as the following simplicial complex: vertices of $\CC(S)$ are
non-trivial, non-peripheral homotopy classes of simple curves, and
$k$--simplices are sets $\{v_0,\ldots,v_k\}$ of distinct vertices with
disjoint representatives. For $k\ge 0$ let $\CC_k(S)$ denote the
$k$--skeleton of $\CC(S)$. 

We define a metric on $\CC(S)$ by making each simplex regular Euclidean
with sidelength 1, and taking the shortest-path  distance.
Distance in $\CC(S)$ will be denoted $d_{\CC(S)}$ or $d_S$. 
Note that $\CC_k(S)$ is quasi-isometric to
$\CC_1(S)$ (with the induced path metrics) for $k\ge 1$.
These
conventions also apply to the non-generic cases below.
\segment{Once-punctured tori and 4--punctured spheres}
If $S$ is $S_{0,4}$ or $S_{1,1}$, we define the
vertices $\CC_0(S)$ as before, but let edges denote pairs
$\{v_0,v_1\}$ which have the minimal possible geometric intersection
number (2 for $S_{0,4}$ and $1$ for $S_{1,1}$).

We remark that in both these cases $\CC(S)$ is isomorphic to the
classical Farey graph in the plane (see eg Hatcher--Thurston
\cite{hatcher-thurston} or  \cite{minsky:taniguchi}).

\segment{Arc complexes}
If $Y$ is a non-annular surface with punctures, let us
also define the larger {\em arc complex} 
$\CC'(Y)$ whose vertices are 
either properly embedded arcs in $Y$, up to homotopy keeping the endpoints in
the punctures, or essential closed curves up to homotopy.
Simplices  correspond to sets of vertices with representatives
that have disjoint interiors.

\segment{Subsurface projections}
By an {\em essential subsurface} of $S$ we shall always mean an open subset $Y$
which is incompressible ($\pi_1$ injects), non-peripheral (not
homotopic into a puncture) and not homeomorphic to $S_{0,3}$.
We usually consider isotopic subsurfaces to be equivalent.

Let $Y$ be a non-annular essential subsurface of $S$.
We can define a ``projection'' 
$$\pi_Y\co \CC'_0(S) \to \CC'(Y)\union\{\emptyset\}$$
as follows (where in the definition of $\CC'(Y)$  we consider each end
of $Y$ to be a puncture):

If $\alpha\in \CC'_0(S)$ has no essential intersections with $Y$
(including the case that $\alpha$ is homotopic to $\boundary Y$) then
define $\pi_Y(\alpha) = \emptyset$.
Otherwise, $\alpha$ intersects 
$Y$ in a collection of disjoint embedded arcs. Keeping only the
essential ones (for example by  taking geodesic
representatives of $\alpha$ and $\boundary Y$ in some hyperbolic
metric), and taking their homotopy classes modulo the ends of $Y$, we
obtain a simplex in $\CC'(Y)$. Let $\pi_Y(\alpha)$ be the 
barycenter of this simplex.

For convenience we also extend the definition of $\pi_Y$
to $\CC'_0(Y)$, where it is the identity map.

We let $d_Y(\alpha,\beta)$ denote
$d_{\CC'(Y)}(\pi_Y(\alpha),\pi_Y(\beta))$, when these projections
are non\-emp\-ty. Similarly $\diam_Y(A)$ 
denotes $\diam_{\CC'(Y)}(\union_{a\in A} \pi_Y(a))$, where
 $A\subset \CC'_0(S)$.

\segment{Annuli} 
Now let $Y$ be an essential annulus in $S$. We will define $\CC'(Y)$
a little differently:
let $\til Y$ denote the annular cover of
$S$ to which $Y$ lifts homeomorphically. We can compactify $\til Y$ to
a closed annulus $\hhat Y$ in a natural way: the universal cover of $S$
can be identified with $\Hyp^2$, which has a natural compactification
as the closed disk, and the covering $\Hyp^2 \to \til Y$ has deck
group $\Z$ which acts with two fixed points on the boundary. The
quotient of the closed disk minus these two points is the closed
annulus $\hhat Y$. Now define $\CC'_0(Y)$ to be the set of all homotopy
classes {\em rel endpoints} of arcs connecting the two boundaries of
$\hhat Y$.
Define $\CC'(Y)$ by putting a simplex between any finite set of arcs that have
representatives with pairwise disjoint interiors.  

We can make the simplices of $\CC'(Y)$ regular Euclidean as
before. Although now $\CC'(Y)$ is infinite-dimensional it is in fact
quasi-isometric to its 1--skeleton $\CC'_1(Y)$, which in turn is
quasi-isometric to the integers $\Z$ with the usual distance
(see \cite{masur-minsky:complex2}). In particular 
it is easy to see that distance in $\CC'_1(Y)$ is determined by
intersection number: If $\alpha$ and $\beta$  are two distinct
vertices in $\CC'_0(Y)$ with
geometric intersection number $i(\alpha,\beta)$ in $\hhat Y$ then
\begin{equation}\label{annulus intersection bounds distance}
d_Y(\alpha,\beta) = 1+i(\alpha,\beta). 
\end{equation}
We can now define $\pi_Y$ in this case as follows: For
 $\alpha\in\CC'_0(S)$,
{\em lift} a representative of  $\alpha$ to
the annular cover $\til Y$. Each component of the lift extends
naturally to an arc in the closed annulus $\hhat Y$, and a finite
number of these connect the two boundary components, and have disjoint
interiors. Thus they determine a simplex of $\CC'(Y)$ and
again we let $\pi_Y(\alpha)$ be the barycenter of this simplex.

\segment{Properties of $\CC(S)$}
We note first that $\CC(S)$ is connected and has infinite diameter in all
cases we consider, and is a $\delta$--hyperbolic metric space for some
$\delta(S)>0$ --- see \cite{masur-minsky:complex1}. 

We will also need an elementary lemma relating $\CC$--distance to
intersection number, analogous to (\ref{annulus intersection bounds
  distance}):

%
%
%

\begin{lemma}{bounded curves close}
Given a surface $Y$ and $D>0$ there exists $D'$ such that, if
$\alpha,\beta$ are vertices of $\CC'(Y)$ and 
$i(\alpha,\beta) \le D$ then
$$
d_Y(\alpha,\beta) \le D'.
$$
\end{lemma}

In fact  a more precise bound can be given --- see
eg \cite{masur-minsky:complex1} and Hempel 
\cite{hempel:complex}. For us it will suffice to note that {\em some}
bound exists. The proof is an easy surgery argument: One can
inductively replace $\beta$ with $\beta'$ which is disjoint from
$\beta$ and intersects $\alpha$ fewer times.

\section{Pleated surfaces}
\label{appendix}
\label{pleated surfaces}

Thurston introduced pleated surfaces in \cite{wpt:notes,wpt:I} as a
powerful tool for studying hyperbolic 3--manifolds. 
In this section we recall the basic definitions of
pleated surfaces and related notions
(for further details see Canary--Epstein--Green
\cite{ceg}), and prove a collection of
results, some of them 
technical and all of them directed towards some extensions and
applications of Thurston's theorems on Uniform Injectivity and Efficiency
of pleated surfaces.

\segment{Definitions}
\label{pleated defs}
See Thurston \cite{wpt:notes}, Canary--Epstein--Green \cite{ceg} and
Penner--Harer  \cite{penner-harer} for 
definitions and examples
of geodesic laminations and measured geodesic laminations
on hyperbolic surfaces. 
 We note here that a geodesic lamination with
respect to one hyperbolic metric on a surface can be isotoped to 
a geodesic lamination with regard to any other hyperbolic metric ---
see Hatcher \cite{hatcher}. If we use the term ``geodesic
lamination'' in the absence of a specific metric we will mean this
isotopy class.

A pleated surface is a map $f\co S\to N$, where $S$ is a surface of
finite type and $N$ is a hyperbolic 3--manifold, together with a
complete, finite-area hyperbolic metric $\sigma$ on $S$, satisfying
the following two conditions. 

First, $f$ is {\em path-isometric} with respect to $\sigma$. That is, 
any $\sigma$--rectifiable path in $S$ is mapped by $f$ to a path of
$N$--length equal to its $\sigma$--length.
With this in mind we call $\sigma$ the
metric {\em induced by $f$}, and note that in fact it is uniquely
determined by $f$.

Second, there is a $\sigma$--geodesic lamination $\lambda$ on $S$ whose
leaves are mapped to geodesics by $f$. In the complement of $\lambda$,
$f$ is totally geodesic.

\paragraph{Double incompressibility}
A map $f\co S\to N$ is {\em incompressible} if $\pi_1(f)$ is injective.
Following Thurston \cite{wpt:I}, we say that 
a map $f\co S\to N$ is {\em doubly incompressible} if in addition
the following conditions hold:
\begin{enumerate}
\item {\em Arcs modulo cusps are mapped injetively:} 

Homotopy classes of maps $(I,\boundary I)\to (S,cusps(S))$
map injectively to 
homotopy classes of maps $(I,\boundary I)\to (N,cusps(N))$.

\item {\em No essential annuli except at parabolics:}

For any cylinder $c\co S^1\times I \to N$ whose boundary 
$\boundary c\co S^1\times\{0,1\} \to N$ factors as 
$f\circ c_0$ where $c_0\co S^1\times\{0,1\} \to S$, if $\pi_1(c)$ is
injective then either the image of $\pi_1(c_0)$ consists of parabolic
elements of $\pi_1(S)$, or $c_0$ extends to a map of $S^1\times I$ 
into $S$.

\item {\em Primitive elements are preserved and no rank--2 cusps in image:}

Each maximal abelian subgroup of $\pi_1(S)$ is mapped to a maximal
abelian subgroup of $\pi_1(N)$.

\end{enumerate}
Here $I=[0,1]$ and $cusps(S)$, $cusps(N)$ are the union of
$\ep_0$--cusp neighborhoods of the parabolic cusps of $S$ and
$N$ respectively. 

In particular, a map which induces an
isomorphism on $\pi_1$ and sends cusps to cusps is doubly incompressible.

\segment{Noded pleated surfaces} Let $S'$ be an essential subsurface
of $S$  whose complement is a disjoint union of simple
curves. Given a homotopy class 
$[f]$ of maps from $S$ to $N$, we say that $g\co S'\to N$ is a {\em noded
  pleated surface} in the class $[f]$ if $g$ is pleated with respect
to a hyperbolic metric on $S'$ (in which the ends are cusps), and $g$ is
homotopic to the restriction to $S'$ of an element of $[f]$.

We note that if $f\co S\to N$ is a doubly-incompressible map and $g\co S'\to
N$ is a noded pleated surface in the homotopy class of $f$, then $g$
is also doubly incompressible.

\segment{Finite-leaved laminations}
In this paper, pleated surfaces and noded pleated surfaces will only
 arise with {\em finite-leaved} laminations,
which always consist of a finite number (possibly zero) of closed
 geodesics and a 
finite number of infinite leaves whose ends either spiral around the
closed geodesics or exit a cusp. 
Thurston observed that, if $f\co S\to N$ is a doubly-incompressible
map and $\lambda$ is a finite leaved lamination all of whose closed
 curves are taken by $f$ to non-parabolic loops, then there is a
 pleated surface homotopic to $f$ mapping $\lambda$ geodesically
 ($\lambda$ is ``realizable'').

Similarly if some subset $C$ of the closed loops are taken by $f$ to
 parabolic loops and $S'=S\setminus C$ then there is a {\em noded} pleated
 surface $g\co S'\to N$ in the homotopy class of $f$.
To simplify notation we will still refer to this noded pleated surface
as a ``pleated surface mapping $\lambda$ geodesically'' (see in
 particular Lemma  \lref{Efficiency of pleated surfaces}).

\segment{Basic properties}
We refer the reader to Thurston \cite{wpt:textbook}, Canary--Epstein--Green
\cite{ceg} or Benedetti--Petronio \cite{benedetti-petronio}  for a discussion 
of the Margulis lemma and the thick--thin decomposition. In what
follows $\ep_0$ will always denote a constant no greater than the
Margulis constant for $\Hyp^3$. We may also assume $\ep_0$ has the
property that the intersection of any simple closed geodesic in a
hyperbolic surface 
with an $\ep_0$--Margulis  tube, is either the core of the tube or 
a union of arcs which connect the two boundaries
(see eg \cite{ceg}).

It is obvious that a $\pi_1$--injective pleated surface $g\co S\to N$ takes the
$\ep_0$--thin part of $S$ to the $\ep_0$--thin part of $N$. In the other
direction we have this observation (see Thurston \cite{wpt:I}):
\begin{lemma}{thin to thin}
Given $\ep_0$ there exists $\ep_1$ (depending only on $\ep_0$ and the
topological type of $S$) such that, if $g\co S\to N$ is a $\pi_1$--injective
pleated surface
then only the $\ep_0$--thin part of the surface can be mapped into the
$\ep_1$--thin part of $N$.
\end{lemma}

\segment{Thurston's uniform injectivity theorem}
One of the most important properties of pleated surfaces is that,
under appropriate topological assumptions, there is some control over
the ways in which they can fold 
in the target manifold. Let $g\co S\to N$ be a pleated surface mapping a
lamination $\lambda$ geodesically, and let $\projmap_g\co \lambda\to \proj(N)$ be
the natural lift of $g|_\lambda$ to the projectivized tangent bundle
$\proj(N)$ of $N$, taking $x\in\lambda$ to $(g(x),g_*(l))$ where $l$ is the
tangent line to the leaf of $\lambda$ through $x$.
In \cite{wpt:II}, Thurston established the
following theorem:
\begin{theorem+}{Uniform Injectivity Theorem}
  Fix a surface $S$ of finite type, and a constant $\bar\ep>0$. 
  Given $\ep>0$ there exists $\delta>0$ such that, given any
  type-preserving doubly incompressible pleated
  surface $g\co S\to N$, with pleating locus $\lambda$ and induced
  metric $\sigma$, if $x,y\in \lambda$ are in the $\bar\ep$--thick part
  of $(S,\sigma)$ then
$$
  d_{\proj(N)}(\projmap_g(x),\projmap_g(y)) \le \delta
\quad \implies\quad
d_\sigma(x,y)  \le \ep.
$$
\end{theorem+}

For more general versions see also
  Thurston \cite{wpt:III}, and Canary \cite{canary:schottky}.)

As a consequence of this theorem, Lemma \ref{Short bridge arcs} below
allows us to compare pleated
surfaces which have a subset of their pleating locus in common.
We first need to address the following somewhat technical point,
arising from the fact that, in the absence of a single fixed
hyperbolic metric, laminations on a surface are usually considered only up
to isotopy.

Let us say that two pleated surfaces $f,g\co S\to N$
are {\em homotopic relative to a common pleating lamination $\mu$} if
$\mu$ is a lamination on $S$ which is mapped geodesically by
$f$, if $f|_\mu = g|_\mu$, and finally if $f$
and $g$ are homotopic by a family of maps that fixes
$\mu$ pointwise. The next lemma tells us that in our setting, if two
pleated surfaces  share a pleating sublamination only up to isotopy, 
then after re-parameterizing the domain the stronger pointwise
condition will hold:

\begin{lemma}{pointwise common pleating}
  Let $f,g\co S\to N$ be homotopic pleated surfaces that are injective
  on $\pi_1$. Suppose that $\mu$ is a sub-lamination of
  the pleating locus of $f$ and $\mu'$ is a sub-lamination of the
  pleating locus of $g$, such that $\mu$ and $\mu'$ are isotopic to
  each other. Then there is a homeomorphism $h\co S\to S$, isotopic to
  the identity, such that $h(\mu) = \mu'$, and $f$ and $g\circ h$ are
  homotopic relative to the common pleating lamination $\mu$.
\end{lemma}

\begin{pf}
By assumption there is some homeomorphism $k\co S\to S$ isotopic to the
identity such that $k(\mu) = \mu'$. Thus possibly replacing $g$ by
$g\circ k$ we may assume $\mu=\mu'$.

Next, after possibly lifting to a suitable cover of $N$ we may assume
that $f$ and $g$ are isomorphisms on $\pi_1$. 
The image of $\mu$ by $f$ is a geodesic lamination $\mu^*$ in $N$ ---
this is also the image of $\mu$ by $g$ since they are homotopic, as below.

Let $H\co S\times I\to N$ be the homotopy
from $f$ to $g$, where $I=[0,1]$ and $H(x,0)=f(x)$, $H(x,1)=g(x)$.
Now, each leaf $l$ of $\mu$ is taken by $f$  to a leaf of
$\mu^*$ in $N$, and similarly by $g$. 
If we lift the homotopy $H$ to a map of
the universal covers $\til H \co  \til S\times I \to \til N$, and let
$\til l$ be a lift of $l$, then we find that the paths $\til H(\til l
\times \{t\})$ for $t\in I$ are all homotopic with bounded
trajectories,  and hence the geodesics 
$\til H(\til l \times \{0\})$ and
$\til H(\til l \times \{1\})$ are the same. 
Projected back down to the manifolds, this means that 
$f$ and $g$ take $l$ to the {\em same} leaf  $l^*$. Furthermore it means
that each trajectory $H(\{x\}\times I)$,
  for $x\in \mu$, is homotopic rel 
endpoints into the image $l^*$ of the leaf $l$ containing $x$. 

Now define a new homotopy $\hat H\co S\times I\to N$, by letting
$\hat H|_{\{x\}\times I}$, for each $x\in S$, be the unique constant
  speed geodesic homotopic to $H_{\{x\}\times I}$ with endpoints
  fixed. This is clearly a 
  continuous map (here we are using the negative curvature of the
  target), and by the previous
  paragraph, when $x$ is in a leaf $l$ of $\mu$, $\hat H(\{x\}\times
  I)$ lies in $l^*$.

For each $t\in I$, let $g_t(\cdot) = \hat H(\cdot,t)$, so that $g_1=g$
and $g_0=f$.
Then $g_t$ maps $\mu$ to $\mu^*$, and in
fact the lift $\projmap_{g_t} \co  \mu \to \proj(N)$ takes $\mu$
homeomorphically to the lift $\proj\mu^*$ in the projective tangent
bundle. 
Let $h_t \co  \mu \to \mu$ be the homeomorphism $\projmap_{g_t}^{-1}
\circ \projmap_f$. Note that $h_0$ is the identity, and
that $h_t$ takes each leaf of $\mu$ to itself. Furthermore
$h_t$ has the property that 
$g_t(h_t(x)) = f(x)$ for all $x\in \mu$.

It remains to extend $h_t$ to a homeomorphism of $S$, and then $g\circ
h_1$ is the desired reparameterization of $g$, with the desired
homotopy being $G(x,t) = g_t(h_t(x))$.

Each complementary region $R$ of $\mu$ is an open hyperbolic
surface with finite area. Taken with the induced path metric,
it has a metric completion $\overline R$, which is a hyperbolic
surface with geodesic
boundary possibly with ``spikes'' --- that is, each end of $R$ is
completed either by a closed geodesic or by a chain of geodesic lines
with each two  successive lines asymptotic along rays. The region
between two such rays where they are less than $\ep$ apart (for some
small $\ep>0$) is called an $\ep$--spike. All these lines naturally immerse in
$S$ as leaves of $\mu$. Each $h_t$ induces a homeomorphism of the boundary of
$\overline R$ which is uniformly continuous (since $h_t$ is continuous on a
compact set), and is easy to extend continuously to an embedding of
each $\ep$--spike into $\overline R$ for suitable $\ep$, which remains
uniformly continuous and continuous in $t$
(for example use horocyclic coordinates
in each spike and extend the map by a linear interpolation of the
boundary values). $\overline R$ minus its spikes is a compact surface and the
standard methods apply to extend $h_t$ to a homeomorphism of
$\overline R$, continuously varying with $t$. These extensions piece together
continuously 
when all the components of $\overline R$ are immersed in $S$, because of the
uniform continuity of each.
\end{pf}

\segment{Applications of uniform injectivity}

If $\alpha$ is a lamination in $S$, a {\em
  bridge arc} for $\alpha$ is an arc in $S$ with endpoints on
$\alpha$, which is not
deformable rel endpoints into $\alpha$. 
A {\em primitive bridge arc} is a bridge arc whose interior is disjoint
from $\alpha$.
If $\sigma$ is a hyperbolic metric on $S$ and $\tau$ is a bridge arc
for $\alpha$,
we let $\ell_\sigma([\tau])$ denote the $\sigma$--length of the
shortest arc homotopic to $\tau$, with endpoints fixed.

\begin{lemma+}{Short bridge arcs}
  Fix a surface $S$ of finite type. Fix a positive $\ep_1\le \ep_0$. 
  Given $\delta_1>0$ there exists
  $\delta_0\in(0,\delta_1)$ such that the following holds. 

  Let $g_0\co S\to N$ and  $g_1\co S\to N$ be 
  type-preserving doubly incompressible pleated surfaces in a hyperbolic
  3--manifold $N$, which are homotopic relative to a common
  pleating lamination $\mu$.
  Let $\sigma_0,\sigma_1$ be the induced metrics
  on $S$. If $\tau$ is a bridge arc for $\mu$, and either
  \begin{enumerate}
  \item $\tau$ is in the $\ep_1$--thick part of $(S,\sigma_0)$, or
  \item $\tau$ is a primitive bridge arc, 
  \end{enumerate}
then
$$
\ell_{\sigma_0}([\tau]) \le \delta_0 \implies
\ell_{\sigma_1}([\tau]) \le \delta_1.
$$
\end{lemma+}

This lemma is a direct consequence of 
Lemma 2.3 in \cite{minsky:3d}, and the discussion of Uniform Injectivity 
there. Let us sketch the proof.

The bound
$\ell_{\sigma_0}([\tau])\le \delta_0$ implies a proportional bound
$d_{\proj(N)}(\projmap_{g_0}(x),\projmap_{g_0}(y))\le C\delta_0$, where $x$ and
$y$ are the endpoints of $\tau$ and $C$ is a universal constant.
This is because, for a uniform $c$, there are segments of 
leaves of $\mu$ centered on $x$ and $y$ of radius at least $\log
{1/\delta_0}-c$, with lifts to the universal cover of
$(S,\sigma_0)$ 
that are bounded Hausdorff distance $\ep_0$ --- these map to 
segments in $N$ with similar bounds since $g_0$ is Lipschitz, 
and hyperbolic trigonometry implies the bound on $d_{\proj(N)}$.

Note also that $d_{\proj(N)}(\projmap_{g_0}(x),\projmap_{g_0}(y))
= d_{\proj(N)}(\projmap_{g_1}(x),\projmap_{g_1}(y))$, since
$g_0|_\mu=g_1|_\mu$ by definition. 

Now consider case (1) first.  
Because $\tau$ is contained in the $\ep_1$--thick part of $(S,\sigma_0)$, by
Lemma \ref{thin to thin} there is some $\ep_2(\ep_1)$ so that $g_0(\tau)$
lies in the $\ep_2$--thick part of $N$. 
Choose $\ep>0$ so that $\ep+\delta_0 < \ep_2$.
Now $x,y$ must be in the $\ep_2$--thick part of $(S,\sigma_1)$, so 
the Uniform Injectivity Theorem applied to  $g_1$, with $\bar\ep = \ep_2$,
gives a choice of $\delta_0$ that guarantees
$d_{\sigma_1}(x,y) \le \ep$. Let
$\tau'$ be an arc in $S$ joining $x$ and $y$ of length $\ep$.
Then the concatenation 
$g_0|_\tau * g_1|_{\tau'}$ has length bounded by $\delta_0 + \ep <
\ep_2$ and is hence 
null-homotopic. By assumption, $g_0$ and $g_1$ are homotopic on $\tau'$
rel endpoints and we conclude that $\tau$ and $\tau'$ are homotopic
rel endpoints. Thus the bound on $\tau'$ yields
the desired bound on $\ell_{\sigma_1}([\tau])$.

In case (2), although $\tau$ can be in the thin part, 
the extra topological restriction that $\tau$ is a
primitive bridge arc is used in Lemma 2.3 of \cite{minsky:3d} to
obtain the desired bound. Essentially, this condition prevents $x$ and
$y$ from being on opposite sides of the core of a
Margulis tube that is very badly
folded in the 3--manifold --- this is the basic example of failure of
Uniform Injectivity. See also Brock \cite{brock:continuity} for a
stronger version of this result.

 A further consequence of Uniform Injectivity allows us to estimate
the lengths of curves in a hyperbolic 3--manifold based on their
representatives in a certain kind of pleated surface. Let $\lambda$ be
a finite-leaved geodesic lamination in a hyperbolic surface $S$ which
is maximal, in the 
sense that its complementary regions are ideal triangles (when $S$ has
cusps there will be leaves that enter the cusps).
If $c$ is any closed geodesic in $S$ we define its {\em alternation
  number} with $\lambda$, $a(\lambda,c)$, following Thurston
\cite{wpt:II}:

\begin{figure}[ht!]
\cl{\relabelbox\small
\epsfbox{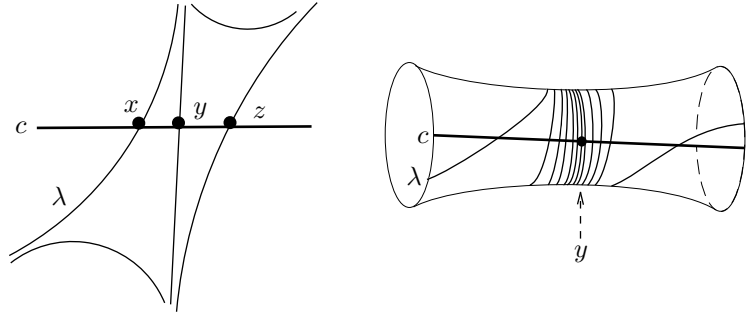}
\relabel{x}{$x$}
\relabel{y}{$y$}
\relabel{y1}{$y$}
\relabel{z}{$z$}
\relabel{c}{$c$}
\relabel{c1}{$c$}
\relabel{l1}{$\lambda$}
\relabela <-2pt, 0pt> {l}{$\lambda$}
\endrelabelbox}
\caption{The two kinds of
  boundary intersections used in defining $a(\lambda,c)$}
\label{boundary intersections}
\end{figure}

If $c$ is a leaf 
of $\lambda$, then  $a(\lambda,c)=0$. If not then $c$ meets $\lambda$
in a countable set of transverse intersection points. 
If $x$ and $y$ are two points of $\lambda\intersect c$ 
bounding an interval on $c$ with no intersection points in its
interior, then $x$ and $y$ are on leaves of $\lambda$ that are legs of
an ideal triangle, and hence are asymptotic, either on one side of $c$
or the other. If $x,y$ and $z$ are successive intersection points with
no intervening intersections then we call $y$ a {\em boundary
  intersection} if the leaves through $x$ and $y$ are asymptotic on
the opposite side of $c$ from the leaves through $y$ and $z$. 
If $y\in\lambda\intersect c$ is an accumulation point of  leaves of
$\lambda$, then since $\lambda$ is finite-leaved the leaf through $y$
is closed, and there are accumulations from each side, by spiraling
leaves. Again we call $y$ a boundary intersection if the spiraling is
in opposite directions on the two sides (see Figure \ref{boundary
  intersections}).

Each boundary intersection is isolated, so there are finitely
many. We let $a(\lambda,c)$ denote the number of boundary intersections.
%

The following theorem will be used in the proof of Theorem B.
The first part (\ref{Efficiency original})
was proved by Thurston in \cite[Theorem 3.3]{wpt:II}.
A generalization was proved by Canary in \cite[Proposition 5.4]{canary:schottky},
to which we refer the reader for a detailed
discussion. The second conclusion (\ref{Efficiency relative}) in our
statement comes from a straightforward extension of Thurston's
argument, which we sketch below.

\begin{theorem+}{Efficiency of pleated surfaces}
Given $S$ and any $\ep>0$, there is a constant $C>0$ such that, if
$g\co S\to N$ is a type-preserving doubly incompressble pleated surface
with induced metric $\sigma$ mapping geodesically a
maximal finite-leaved  lamination $\lambda$, and each closed leaf of
$\lambda$ has image length at least $\ep$, then for any measured
geodesic lamination $\gamma$ in $S$,
\begin{equation}\label{Efficiency original}
\ell_N(g(\gamma)^*) \le \ell_\sigma(\gamma) \le \ell_N(g(\gamma)^*)
+ Ca(\lambda,\gamma).
\end{equation}
If we remove the length condition on the closed leaves of $\lambda$ and let
$R_\lambda$ be the complement in $S$
of the $\ep_0$--thin parts (in the metric $\sigma$)
whose cores are closed leaves of
$\lambda$, we have
\begin{equation}\label{Efficiency relative}
\ell_\sigma(\gamma\intersect R_\lambda) \le
\ell_N(g(\gamma)^*) + Ca(\lambda,\gamma).
\end{equation}
This estimate applies also if some  closed leaves of $\lambda$ have zero
length, in which case $g$ is a noded pleated surface defined in the
complement of those leaves.
\end{theorem+}
Here $g(\gamma)^*$ denotes the geodesic representative of $g(\gamma)$
in $N$, if it exists. For a measured lamination $\gamma$ this means
the (unique) image of $\gamma$ by a pleated surface homotopic to $g$
that maps $\gamma$ geodesically. The length $\ell_N(g(\gamma)^*)$ is then
well-defined, and if no geodesic representative exists ($\gamma$ is
mapped to a parabolic element, or is an ending lamination) then we
define $\ell_N(g(\gamma)^*)=0$. This quantity is continuous as  a
function of the measured lamination $\gamma$ (see Thurston \cite{wpt:II} and
also Brock \cite{brock:continuity}).

\begin{pf}(Sketch)\qua
Let us first recall Thurston's 
original argument for (\ref{Efficiency original}).
We may assume that $\gamma$ is a simple closed curve --- the general case can
be obtained by taking limits, since multiples of simple closed curves
are dense in the measured lamination space and the length function is
continuous. 
The first step is to construct a polygonal curve $\mu$ on $S$, homotopic to
$\gamma$, which consists of $a(\lambda,\gamma)$ segments on leaves of
$\lambda$ (which meet $\gamma$ at the boundary intersections),
connected by $a(\lambda,\gamma)$ ``jumps'' of bounded length. 

This is best seen by lifting $\gamma$ to a line $\til\gamma$ in $\Hyp^2$ and
considering the chain of geodesics $\{g_i\}_{i=-\infty}^\infty$
in the lift of $\lambda$ that
cross it in boundary intersections. Thus $g_i$ and $g_{i+1}$ are
asymptotic on alternating sides of $\til\gamma$. 
Let $X$ be a 1--neighborhood of $\til\gamma$. A path $\til\mu$ is
constructed as a chain of segments, alternately subsegments of
$g_i\intersect X$ and segments in $X$ of uniformly bounded length joining 
$g_i\intersect X$ to $g_{i+1}\intersect X$. See Figure \ref{polygonal
  approx} for an example, and Thurston or Canary for the exact
construction. 

\begin{figure}[ht!]
\cl{\relabelbox\small
\epsfbox{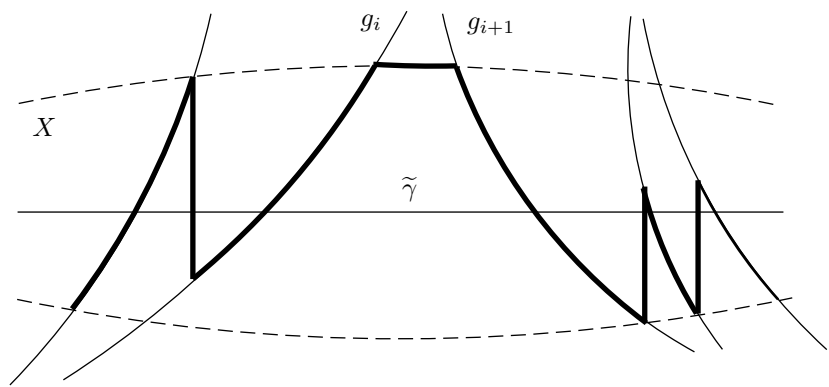}
\relabel{X}{$X$}
\relabel{g}{$\widetilde\gamma$}
\relabel{gi}{$g_i$}
\relabel{g1}{$g_{i+1}$}
\endrelabelbox}
\caption{Thurston's polygonal
  approximation to a geodesic using the leaves of a maximal
  finite-leaved lamination. The boundary of $X$ is dotted; $\widetilde\mu$
  is thickened.}
\label{polygonal approx}
\end{figure}

We project $\til\mu$ to $S$ to obtain a closed polygonal curve $\mu$ with
$2a(\lambda,\gamma)$ segments. 
This curve has the property that 
$\ell_\sigma(\mu) \le \ell_\sigma(\gamma) + C_0a(\lambda,\gamma)$.
(In the rest of the argument each
$C_i$ denotes some constant independent of anything but $\ep$ and the
topological type of $S$.) 

Furthermore, Thurston points out that we can adjust $\mu$, by moving
the points where the bounded jumps occur, so that 
those jumps never occur inside the $\ep_0$ Margulis
tubes whose cores have length $\ep$ or less and are not leaves of
$\lambda$ (in Thurston's setting  all the Margulis tubes with core
$\ep$ or less are not leaves of $\lambda$).

The image of $\mu$ by $g$ consists
of the images of arcs along $\lambda$, which are still geodesic, and
of the jumps, which still have bounded lengths. 
After straightening out
the images of the bounded jumps we obtain a polygonal curve $\nu$ in $N$, and
we have the estimates
\begin{equation}\label{estimate mu nu}
\ell_\sigma(\gamma) \le 
\ell_\sigma(\mu) \le
\ell_N(\nu) + C_1 a(\lambda,\gamma).
\end{equation}
Thurston next forms a ``pleated annulus'' connecting $\nu$ to the
geodesic representative $g(\gamma)^*$. This is a pleated map of a
hyperbolic annulus $A$ with one geodesic boundary component that maps
to $g(\gamma)^*$, and one polygonal boundary component that maps to $\nu$.
The Gauss--Bonnet theorem bounds the area of
$A$ in terms of the number of corners of $\nu$, which is
$2a(\lambda,\gamma)$. Let us partition $\nu$ into several parts: 
Fixing some $\ep'>0$, 
$\nu_0$ is the union of segments that admit collar neighborhoods of
width $\ep'$ in $A$, together with all the bounded length jumps.
The length of $\nu_0$ is bounded by $C_2a(\lambda,\gamma)/\ep'$ because of
the area bound on $A$. 
Any component of  $\nu\setminus\nu_0$ must run within distance $\ep'$
(in $A$) of some 
other segment of $\boundary A$. Let $\nu_1$ denote those segments 
that are distance $\ep'$ from the boundary component mapping to
$g(\gamma)^*$ --- the total length of these is at most 
$\ell_N(g(\gamma)^*)$ plus a constant times the number of segments in
$\nu_1$, which again is at most $C_3a(\lambda,\gamma)$. 

Let $\nu_2$ denote the rest of $\nu$ --- any component of $\nu_2$ 
is within $\ep'$ in $A$ of another segment of $\nu$.
We then bound the length of $\nu_2$ by $C_4 a(\lambda,\gamma)$,
using the Uniform Injectivity Theorem: If two long segments of $\nu$ 
along $\lambda$ are nearly parallel in $A$ and hence in $N$, then
there are three 
possibilities: 

If the segments are in the $\ep_0$--thick part, then the Uniform
Injectivity Theorem implies they are nearly parallel in $S$ as well ---
In fact (see Lemma \ref{Short bridge arcs} for such an argument), with
$\ep'$ chosen sufficiently small compared to $\ep_0$, 
the short arc
connecting them in $A$ is homotopic to a short arc in $S$, so
the curve $\mu$ can be shortened significantly by a homotopy in $S$,
reducing its length by at least the sum of the lengths of the two parallel
segments. However, since
$\ell_\sigma(\gamma) \le
\ell_\sigma(\nu) \le \ell_\sigma(\gamma) + C_5 a(\lambda,\gamma)$,
this is a contradiction unless
the segments involved have length bounded by $C_6a(\lambda,\gamma)$.

If one of the segments is 
in a Margulis tube of $S$ whose core is not a component of $\lambda$, then
both must be in the same Margulis tube, since $g$ takes thin parts to
thin parts, bijectively by the condition on $\pi_1$.
Thus
one can follow both segments until they exit the tube (since all jumps
happen outside this tube) 
and where they exit one can again
use Uniform Injectivity to show the leaves are close in $S$.
Again we obtain too much shortening of $\nu$ if the segments are too
long.

The last possibility,
that the segments are in a Margulis tube whose core has length shorter
than $\ep$ and is a component of
$\lambda$, is disallowed by Thurston's hypothesis. Thus he obtains
$$
\ell_N(\nu) \le \ell_N(g(\gamma)^*) + C_7 a(\lambda,\gamma)
$$
from which (\ref{Efficiency original}) follows.

Now in order to obtain (\ref{Efficiency relative}), we allow the
possibility of Margulis tubes of $\sigma$ whose cores are closed
curves in $\lambda$, and let $R_\lambda$ be the complement of these 
tubes. Let us first assume that no closed leaf of $\lambda$ maps to a
parabolic curve in $N$, and return to this case (the noded case) at
the end.

Thurston's argument applies to all parts of $\nu$ in $R_\lambda$,
and hence immediately yields
\begin{equation}\label{estimate rel nu gamma}
\ell_N(\nu\intersect R_\lambda) \le \ell_N(g(\gamma)^*) + C_8
a(\lambda,\gamma).
\end{equation}
We also have 
\begin{equation}\label{estimate rel mu nu}
\ell_\sigma(\mu\intersect R_\lambda) \le \ell_N(\nu\intersect
R_\lambda) + C_1a(\lambda,\gamma)
\end{equation}
by the same argument as for (\ref{estimate mu nu}).
It remains to establish a connection between 
$\ell_\sigma(\mu\intersect R_\lambda)$ and  
$\ell_\sigma(\gamma\intersect R_\lambda)$.
To do this, we note that the configuration of $\lambda$ in a Margulis
tube outside $R_\lambda$ is as follows: There is the core curve $c$ of the
tube, which is a leaf of $\lambda$, and the other leaves of $\lambda$ that
meet the tube must enter from one side or the other and spiral around
$c$. An arc of $\gamma$ that runs through this tube 
has one boundary intersection each with two leaves of $\lambda$,
on opposite sides of the core $c$, and possibly
one boundary
intersection with the $c$ itself (if the spiraling on the two sides of
$c$ is in opposite directions).
Thus the intersection of $\mu$ with this
Margulis tube consists of at most three arcs of $\lambda$, each of 
which remains in a bounded neighborhood of
$\til\gamma$ in the lift to $\Hyp^2$, 
and at most two bounded jumps.
It follows that the arc of  $\mu$ intersected with this tube has
length at most a
constant plus the length of the corresponding arc of $\gamma$.
We conclude that
$$
\ell_\sigma(\gamma\intersect R_\lambda) \le
\ell_\sigma(\mu\intersect R_\lambda) + C_9 a(\lambda,\gamma).
$$
This together with (\ref{estimate rel nu gamma}) and 
(\ref{estimate rel mu nu})
gives the desired inequality.

Finally, we consider the noded case, where some closed leaves
$\{c_1,\ldots c_k\}$ of
$\lambda$ are mapped to parabolics in $N$. Then we let $S'$ be the complement
of these leaves and $g\co S'\to N$ a noded pleated surface mapping
geodesically all leaves of $\lambda$ except the $c_i$. The $g$--image
of an end of a leaf of $\mu$ 
that spirals around $c_i$ must then terminate in the corresponding
cusp of $N$, and 
two such ends of leaves have asymptotic images. If $\gamma$ does not
intersect any of the $c_i$ the argument can be repeated as before, but
in general $\gamma\intersect S'$ may be a union of arcs, each
represented by an infinite geodesic in $g(S')$ with its ends in the
cusps. The construction of $\mu$ in $S$ is the same as before, but
where $\mu$ crosses one of the $c_i$ it
leaves the domain of $g$. 
For each such crossing, the image $g(\mu\intersect S')$ traverses two
leaves of $\lambda$ that are 
asymptotic into the corresponding cusp. To build $\nu$, we simply
extend these leaves far enough into the cusp that we can jump between
them with a bounded arc in the correct homotopy class. 
As before Thurston's argument applies to all parts of $\nu$ in
$R_\lambda$, which excludes the segments in the parabolic cusps, and
the length of $\gamma\intersect R_\lambda$ and $\mu\intersect
R_\lambda$  again differ by a bounded multiple of $a(\lambda,\gamma)$.
\end{pf}

\section{The proof of Theorem B}
\label{proof of B}

Recall the statement of our second main theorem:
\bfheading{Theorem B}{\sl
Given a surface $S$, $\ep>0$ and $L>0$, there exists $K$ 
so that, if $\rho\co \pi_1(S)\to \PSL 2(\C)$ is a Kleinian surface group
and $Y$ is a proper essential subsurface of $S$, then 
$$
\diam_Y(\CC_0(\rho,L)) \ge K
\quad\implies\quad
\ell_\rho(\boundary Y) \le \ep.
$$
}

In fact, we will prove the equivalent conclusion
$$
\ell_\rho(\boundary Y) >  \ep
\quad\implies\quad
\diam_Y(\CC_0(\rho,L)) < K.
$$
We first give the proof in the case where $Y$ is not an annulus.
The annular case is similar, and we describe it 
at the end of the section.

\segment{The proof in the non-annular case}
\label{non-annular}
The proof reduces to two main lemmas. 
First, given a hyperbolic surface $(Z,\sigma)$ with geodesic boundary,
we define a  ``minimal proper arc''  to be an embedded arc
$\tau$ with endpoints on $\boundary Z$, not homotopic into $\boundary Z$,
which is minimal in $\sigma$--length
among all such arcs.

\begin{lemma}{Y distance to tau}
  For $L>0$ there exists $D>0$, depending only on $L$ and the topology
  of $S$, such that the following holds:

Given a Kleinian surface group $\rho\co \pi_1(S)\to\PSL 2(\C)$,
a non-annular proper essential subsurface $Y\subset S$,  and
any simple closed curve $\gamma$ in $S$ intersecting $Y$ essentially,
  such that $\ell_\rho(\gamma) \le L$,  there is a pleated
  surface $g_\gamma$ in the homotopy class of $\rho$  mapping
  $\boundary Y$ geodesically,  
  with induced metric $\sigma_\gamma$, such that for any minimal
  proper arc $\tau$ in $(Y,\sigma_\gamma)$ we have
\begin{equation}
\label{bound gamma tau}
d_Y(\gamma,\tau)\le D.
\end{equation}
\end{lemma}
We remark that in this lemma we do not use the
assumption on $\ell_\rho(\boundary Y)$.

The second lemma will show that,  over all pleated surfaces
$g\co S\to N$ mapping  
$\boundary Y$ geodesically the set of minimal proper arcs in $Y$ with
respect to the induced metrics has uniformly bounded
diameter in $\CC'(Y)$. More precisely:

\begin{lemma}{pleating on common boundary}
For any $\ep>0$ there 
exists $D$ (depending on $\ep$ and the topological type of $S$) such
that the following holds. 

Let $\rho\co \pi_1(S)\to\PSL 2(\C)$ be a Kleinian surface group, and
let $Y$ be a non-annular proper essential subsurface of $S$ and
$g_0,g_1$ a pair of 
pleated surfaces in the homotopy class $[\rho]$
mapping $\boundary Y$ to
geodesics. Let $\sigma_0$ and $\sigma_1$ be the induced metrics on $S$.
Suppose that $\ell_\rho(\boundary Y) \ge \ep>0$. 
If $\tau_0$ and $\tau_1$ are 
minimal proper arcs
in $(Y,\sigma_0)$ and $(Y,\sigma_1)$, respectively, then
$$
d_Y(\tau_0,\tau_1) \le D.
$$
\end{lemma}

Thus, assuming $\ell_\rho(\boundary Y) \ge \ep$, 
for any two curves $\gamma_0$, $\gamma_1$ in $\CC_0(\rho,L)$
whose projections $\pi_Y(\gamma_i)$ are non-empty, we
apply Lemma \ref{Y distance to tau} to obtain two pleated surfaces
mapping $\boundary Y$ 
geodesically and minimal proper arcs $\tau_0$, $\tau_1$ in $Y$ with respect
to the two induced metrics, with a bound on $d_Y(\gamma_i,\tau_i)$
from (\ref{bound gamma tau}). Lemma \ref{pleating on common boundary}
then implies the final bound on $d_Y(\gamma_0,\gamma_1)$. This
completes the proof modulo the two lemmas.

Before giving the proofs let us recall the following geometric fact:
Let $Z$ be a complete hyperbolic or simplicial hyperbolic surface with
boundary curved outward (curvature vector pointing out of the surface)
--- for example a surface with geodesic boundary or a cusped surface
minus a horoball neighborhood of the cusp.
Then the length $r$ of  a minimal proper
arc in $Z$ satisfies  
\begin{equation}
\label{collar bound}
2\pi|\chi(Z)| \ge \ell(\boundary Z) \sinh (r/2),
\end{equation}
by the Gauss--Bonnet theorem and an elementary formula for the area of
an embedded collar around $\boundary Z$.

\begin{pf*}{Proof of Lemma \ref{Y distance to tau}}\,
Let $\gamma$ denote any simple closed curve
with $\ell_\rho(\gamma) \le L$, which has
non-trivial intersection with $Y$. If $\gamma\subset Y$ then we can
simultaneously pleat along $\boundary Y$ and $\gamma$ yielding a
pleated surface $g_\gamma$ in which $\gamma$ has length bound
$L$. 

Let $R_\gamma$ denote the complement in $S$ of the $\ep_0$--Margulis tubes
(with respect to the induced metric of $g_\gamma$)
whose cores are components of $\boundary Y$.

Note also that we allow  
$g_\gamma$ to be a {\em noded} pleated surface, in case any components
of $\boundary Y \union \gamma$ are parabolic in $N$ (see Section
\ref{pleated surfaces}). In this case $g_\gamma$ is defined on the complement
$S'$ of the parabolic curves, and $\tau$ is a geodesic 
arc in $Y$ with endpoints in the cusps and whose length within
$R_\gamma$ is minimal. 

Inequality (\ref{collar bound}) applied to $R_\gamma\intersect Y$ gives a
uniform upper bound on the  length of $\tau\intersect R_\gamma$. This
gives us a  uniform 
upper bound on the intersection number $i(\gamma,\tau)$, and hence
on $d_Y(\gamma,\tau)$ by Lemma \ref{bounded curves close}.

From now on assume that $\gamma$ has nontrivial intersection with
$\boundary Y$. 

\begin{figure}[ht!]
\cl{\relabelbox\small
\epsfbox{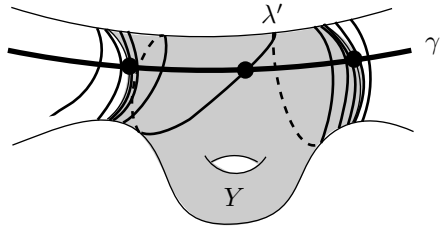}
\relabel{Y}{$Y$}
\relabel{l}{$\lambda'$}
\relabel{g}{$\gamma$}
\endrelabelbox}
\caption{An intersection of $\gamma$ with $Y$ and
  the resulting leaves of the spun lamination $\lambda'$. $Y$ is
  shaded and the resulting boundary intersections are indicated by bullets.
}
\label{spin gamma}
\end{figure}

Let ${\DD}$ denote the mapping class that performs one positive Dehn
twist on each component of $\boundary Y$. The sequence of curves
$\DD^n (\gamma)$ converge, as $n\to \infty$, to a finite-leaved lamination
$\lambda'$ whose non-compact leaves spiral around $\boundary Y$,
in such a way that the spiraling on opposite sides of any component
is in opposite directions. Its closed leaves are exactly 
$\boundary Y$. (This construction is often called ``spinning''
$\gamma$ around $\boundary Y$ --- see Figure \ref{spin gamma} for an example).

After adding a finite number of other infinite leaves spiraling
around $\boundary Y$,
we can obtain a {\em maximal} finite-leaved lamination $\lambda$
containing $\lambda'$. We immediately observe that
\begin{equation}
  \label{a is i}
  a(\lambda,\gamma) = a(\lambda',\gamma) = 2i(\boundary Y,\gamma)
\end{equation}
since $\gamma$ can be represented by a chain of segments of
infinite leaves of $\lambda$, each asymptotic to a component of
$\boundary Y$ at each 
end, with successive segments joined by paths along $\boundary Y$. After
perturbation to a curve transverse to $\lambda$ we obtain a curve
which has one boundary
intersection for each intersection of $\gamma$ with $\boundary Y$, and
one boundary intersection with each of the leaves in the chain.

Let $g_\lambda$ be the pleated surface mapping
$\lambda$ geodesically, 
and let $\sigma_\lambda$ be its induced metric.
Again if some components of $\boundary Y$ are parabolic in $N$
then $g_\lambda$ will be a noded pleated surface defined on the
complement $S'$ of these components, and $\sigma_\lambda$ will have
cusps corresponding to the ends of $S'$. Otherwise let $S'=S$.

Again let
$R_\lambda$
be the complement in $(S',\sigma_\lambda)$  of the $\ep_0$--Margulis
tubes whose cores are components of $\lambda$ (hence of $\boundary Y$).
Realize $\gamma$ by its geodesic representative in
$(S',\sigma_\lambda)$ (which we continue to call $\gamma$; in the
noded case $\gamma$ is replaced by the geodesic representative rel
cusps of $\gamma\intersect S'$).
Theorem \ref{Efficiency of pleated surfaces} part (\ref{Efficiency
  relative})  then implies
\begin{equation}
  \label{consequence of efficiency}
\ell_{\sigma_\lambda}(\gamma\intersect R_\lambda) \le \ell_\rho (\gamma)  + C 
a(\lambda,\gamma).
\end{equation}
Using (\ref{a is i}) and the fact that
$\ell_\rho(\gamma) \le L$, we conclude that
\begin{equation}
  \label{intersection upper bound}
\ell_{\sigma_\lambda}(\gamma\intersect R_\lambda) \le L +
2Ci(\gamma,\boundary Y).   
\end{equation}
We remark that, in case all components of $\boundary Y$ have length at
least $\ep_0$, we may use Thurston's original version of Theorem
\ref{Efficiency of pleated surfaces} (\ref{Efficiency original}), or
equivalently set $R_\lambda = S$.

Let $\gamma_1,\ldots,\gamma_m$ be the component arcs of
$\gamma\intersect (R_\lambda\setminus \boundary Y)$,
where
$m=i(\gamma,\boundary Y)$.  At least half of these arcs are contained
in $Y\intersect R_\lambda$. A bit of algebra shows that, 
choosing $K = 4C + 2L$, less than
$m/2$ of the arcs can have length greater than $K$. Thus we are
guaranteed at least one arc $\gamma_i$ in $Y\intersect R_\lambda$
with $\sigma_\lambda$--length bounded by $K$. 

The length of $\tau\intersect R_\lambda$ is bounded above as before by
(\ref{collar bound}), and hence again we obtain a bound on the
intersection number $i(\gamma_i,\tau)$.
Lemma \ref{bounded curves close} then bounds $d_Y(\gamma_i,\tau)$,
which is within 1 of $d_Y(\gamma,\tau)$ since $\gamma_i$ is in the
same simplex of $\CC'(Y)$ as $\pi_Y(\gamma)$.
\end{pf*}

\begin{pf*}{Proof of Lemma \ref{pleating on common boundary}}\,
We need to recall an additional geometric fact:
There is a constant $\delta_1>0$ for which, if $Z$ is any hyperbolic
surface with geodesic boundary and  $\tau$ and $\tau'$ 
are essential properly embedded arcs in $Z$ whose lengths are
at most $\delta_1$, then $\tau$ and $\tau'$ are either homotopic keeping
endpoints in $\boundary Z$, or they are disjoint. (This follows
directly from the 
fact that boundary components of a hyperbolic surface which are close
at one point must be nearly parallel for long stretches.)

Given this constant $\delta_1$, let $\delta_0\le \delta_1$ be the
constant given by Lemma 
\lref{Short bridge arcs}. 
Note that, by Lemma \ref{pointwise common pleating}, we may assume,
possibly precomposing $g_1$ by a domain homeomorphism isotopic to the
identity, that $g_0$  and $g_1$ are homotopic relative to the common
pleating lamination $\boundary Y$. This allows us to apply Lemma
\ref{Short bridge arcs}.

There are now two cases. 

\bfheading{Case a:} Suppose $\ell_{\sigma_0}(\tau_0)\le \delta_0$. Since
$\tau_0$ is a primitive bridge arc for $\boundary Y$, Lemma \ref{Short
  bridge arcs} part (2) guarantees that $\ell_{\sigma_1}(\tau_0) \le
\delta_1$ as well. Then $\ell_{\sigma_1}(\tau_1) \le \delta_1$ since
$\tau_1$ is minimal in the metric $\sigma_1$, and we may then conclude,
by choice of $\delta_1$, that the arcs $\tau_0$ and
$\tau_1$ are either homotopic or disjoint, and hence
$d_{\CC'(Y)}(\tau_0,\tau_1)\le 1$.

\bfheading{Case b:} Suppose $\ell_{\sigma_0}(\tau_0) > \delta_0$. 
Then (\ref{collar bound}) implies that there is some
$A(\delta_0)$ such that
$\ell_N(\boundary Y) = \ell_{\sigma_0}(\boundary Y) \le A$. 

The rest of the proof resembles that of the Connectivity Lemma
8.1 in \cite{minsky:torus}:

Let $g_t, t\in[0,1]$ be a continuous family  of maps connecting $g_0$
to $g_1$ so 
that for each $t\in(0,1)$, $g_t$ is a simplicial hyperbolic map, and maps
$\boundary Y$ geodesically. Recall that this means that $g_t$ may be
noded on parabolic components of $\boundary Y$, but not on all of them
since we assume $\ell(\boundary Y) \ge \ep > 0$.

The existence of such a family follows from the techniques of Thurston
in \cite{wpt:notes} and Canary \cite{canary:covering}, for example. The 
pleated maps $g_0,g_1$ may be approximated by simplicial hyperbolic
surfaces, in which lamination leaves that spiral an infinite number of times
along the geodesics $\boundary Y$ are replaced by triangulation edges
that terminate on $\boundary Y$. Any two such triangulations may be
connected by a sequence of elementary moves in which the edges on
$\boundary Y$ are fixed and the other edges are replaced one-by-one,
by a theorem of Hatcher \cite{hatcher:triangulations}.
Each such triangulation gives rise to a
simplicial hyperbolic surface in which the components of $\boundary Y$
are mapped to their geodesic representatives (or the surface is noded
on the ones that are parabolic), and
Canary shows in \cite{canary:covering} that each elementary move
between two such surfaces may be realized by a continuous family of 
simplicial hyperbolic surfaces.

The induced metrics $\sigma_t$ vary
continuously, in the sense that for any fixed
homotopy class $\beta$ of curves or arcs rel boundary,
$\ell_{\sigma_t}(\beta)$ is continuous in $t$.

Now given any essential properly embedded arc $\tau$ in $Y$,
let $E_\tau\subset [0,1]$ denote the set of $t$--values for which
$\tau$ is (homotopic rel $\boundary Y$ to)
a minimal proper arc with respect to $\sigma_t$.
Since we are assuming $\ell_\rho(\boundary Y) \ge \ep$, (\ref{collar
bound})  gives an
upper bound $\transL=\transL(\ep)$ on $\ell_{\sigma_s}([\tau])$ for
$s\in E_\tau$.
Continuity of
the metrics implies that $E_\tau$ is closed, and
clearly  the family $\{E_\tau\}$ covers $[0,1]$.

We now observe two facts. 

First,
if $E_\tau \intersect E_{\tau'}\ne \emptyset$ then
$\tau$ and $\tau'$ are simultaneously (homotopic to) shortest arcs in
the same metric 
$\sigma_s$, $s\in E_\tau \intersect E_{\tau'}$. We may assume,
possibly after homotopy, that 
$\tau$ and $\tau'$ are shortest-length representatives of their classes.
We claim that $i(\tau,\tau') \le 1$. If not, let $x$ and $y$ be two
intersection points of $\tau$ and $\tau'$, cutting each into three
successive arcs $\tau_1,\tau_2,\tau_3$ and 
$\tau'_1,\tau'_2,\tau'_3$, respectively. The concatenations
$\xi_1=\tau_1*\tau'_1$ and $\xi_3=\tau_3*\tau'_3$ are
arcs with endpoints in $\boundary Y$, 
and are not homotopic into the boundary: if, say, $\xi_1$ were,
then since it meets $\boundary Y$ at right angles by the minimality of
$\tau$ and $\tau'$, we would obtain a disk that violates the
Gauss--Bonnet theorem. (Actually if the metric $\sigma_t$ is singular 
exactly at a point where $\xi_1$ meets $\boundary Y$ then there is
an angle of {\em at least} $\pi/2$ on each side, and the same argument holds.)
Let $a_i$ and $a'_i$ be the lengths of $\tau_i$ and
$\tau'_i$, respectively. Let $L=a_1+a_2+a_3 = a'_1+a'_2+a'_3$.
Then since $a_2,a'_2>0$ we have $a_1 + a_3 + a'_1 + a'_3 < 2L$, so 
either $a_1 + a'_1<L$ or
$a_3+a'_3<L$. This means that either $\xi_1$ or $\xi_3$ is strictly
shorter than $\tau$ and $\tau'$, contradicting their minimality.

Thus $i(\tau,\tau') \le 1$, and 
this easily gives $d_{\CC'(Y)}(\tau,\tau') \le 2$.

The second fact is that, since $\ell_\rho(\boundary Y) \le A$, 
all the arcs $\tau$ for
which $E_\tau\ne \emptyset$, together with the boundary components of
$Y$ on which their endpoints lie, can be realized in the quotient
hyperbolic 3--manifold $N_\rho$ by a 1--complex with at most $k$
components (where $k$ is the number of
components of $\boundary Y$), each of which has 
diameter at most $A + \transL(k-1)$.
Each $\tau$ together with one or two
segments on $\boundary Y$ gives rise to a loop in this 1--complex
of length at most $A+2\transL$, and the loops for homotopically
distinct $\tau$'s are homotopically distinct.
A standard  application of the Margulis lemma gives an upper bound
$M=M(R)$ for any non-elementary Kleinian group $\Gamma$ 
on the number of distinct elements of $\Gamma$ that can translate any
one point a distance $R$ or less. This gives 
a bound $M'=M'(A,\transL)$ on the
number of (homotopically distinct) arcs $\tau$ with
$E_\tau\ne\emptyset$. 

Now consider the 
graph $T$ whose vertices are those homotopy classes of $\tau$ with
$E_\tau\ne 
\emptyset$, and whose edges are those $([\tau],[\tau'])$ for which
$E_\tau\intersect E_\tau'\ne \emptyset$. The $\CC'(Y)$--distance
between two endpoints 
of an edge is at most $2$ by the first observation, and the
number of vertices is at most $M'$, by the second. 
The fact that the $E_\tau$
cover $[0,1]$ and are closed and finite in number means that $T$ is
connected. 
The diameter of the set of vertices of $T$ is therefore bounded
by $2 (M'-1)$.

For our original arcs $\tau_0$ and $\tau_1$ we have $i\in E_{\tau_i}$ 
and therefore $d_Y(\tau_0,\tau_1) \le 2(M'-1)$.

This concludes the proof of Lemma \ref{pleating on common boundary}
and hence of Theorem B in the non-annular case. 
\end{pf*}

\segment{The annular case} 
\label{annular}
The proof again reduces to two lemmas analogous to Lemmas \ref{Y
  distance to tau} and \ref{pleating on common boundary}, and their
  proofs are similar. We introduce a new bit of notation and discuss
  the differences in the proofs:

Let $Y$ be an annulus and $\alpha$ its core
curve.  If $\alpha$ is geodesic in some hyperbolic surface $S$
we can consider a {\em minimal curve crossing $\alpha$} to be a curve
$\beta$ constructed as follows: Pick one side of $\alpha$ in $S$ and 
let $\tau$ be a minimal length primitive bridge
arc for $\alpha$ that is incident to it on this side. Let $\tau'$ be
a minimal length primitive bridge
arc for $\alpha$ that is incident to it on the other side. If $\alpha$
is non-separating and $\tau$ meets it on both sides then we let
$\tau'=\tau$. Let $\beta$
be a minimal length shortest simple curve  that can be represented as a
concatenation of $\tau$, $\tau'$ (if they are different),  and arcs on
$\alpha$. In particular $\beta$ crosses $\alpha$ once or twice. 

The analogue of Lemma \ref{Y distance to tau} is now the following:
\begin{lemma}{Y distance annulus}
Given $\ep>0$ and $S$ there exists $M>0$ such that, for any
Kleinian surface group $\rho\co \pi_1(S)\to\PSL 2(\C)$, the following
holds:

Let $Y$ be an essential annulus with core $\alpha$ such that
  $\ell_\rho(\alpha)\ge \ep$, and let 
  $\gamma\in\CC_0(\rho,L)$ intersect $\alpha$ essentially. Then
there exists a pleated surface
$g_\gamma$ in
whose induced metric $\sigma_\gamma$ a minimal curve $\beta$ crossing
$\alpha$ satisfies 
\begin{equation}
  \label{gamma beta bound}
  d_Y(\gamma,\beta) \le M
\end{equation}
\end{lemma}


\begin{pf}
The construction of $g_\gamma$ is the same as before (starting by spinning
$\gamma$ around $\alpha$ to obtain a finite-leaved lamination
$\lambda'$), but now we find 
that if $\alpha$ is very short the conclusion is false: $\gamma$ could
wind arbitrarily inside the Margulis tube of $\alpha$ without
violating Theorem \lref{Efficiency of pleated surfaces}. Note that in
the non-annular case $\gamma$ could also have wound around $\boundary Y$ but
this had no effect on the distance $d_Y$. At any rate we now must use
the hypothesis of Theorem B that $\ell_\rho(\alpha) \ge \ep$.

\begin{figure}[ht!]
\cl{\relabelbox\small
\epsfbox{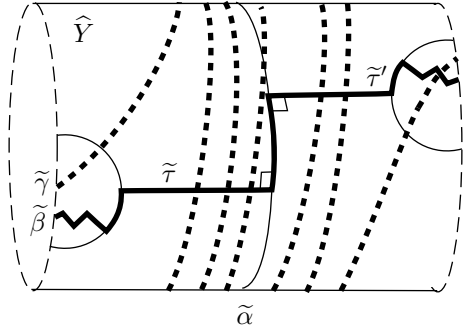}
\relabel{Y}{$\hhat Y$}
\relabela <-2pt, 0pt> {g}{$\til\gamma$}
\relabela <-2pt, 0pt> {b}{$\til\beta$}
\relabela <-2pt, 0pt> {a}{$\til\alpha$}
\relabela <-2pt, 0pt> {t}{$\til\tau$}
\relabela <0pt, 2pt> {t1}{$\til\tau'$}
\endrelabelbox}
\caption{A component $\til\gamma$ of
  the lift of $\gamma$ to the annulus $\hhat Y$ is indicated by the
  broken line. A 
  component $\til\beta$ of the lift of $\beta$ is indicated in heavy
  lines, using a representative that traces along the lifts of
  $\alpha,$ $\tau$ and $\tau'$. Note that $\til\tau$ and $\til\tau'$
  terminate on translates of $\til\alpha$, pictured as semicircles.}
\label{twist beta gamma}
\end{figure}

In this case, part (\ref{Efficiency original}) of Theorem
\lref{Efficiency of pleated surfaces}, together with 
(\ref{a is i}), give us
\begin{equation}
  \label{annulus efficiency}
  \ell_{\sigma_\gamma}(\gamma) \le L + 2Ci(\alpha,\gamma).
\end{equation}
We next relate $d_Y(\beta,\gamma)$ to the number of
times $\gamma$ intersects the arcs $\tau $ and $\tau'$: lift $\alpha$
to the core $\til\alpha$ of the 
annular cover $\hhat Y$ corresponding to $Y$, and consider a
component $\til\gamma$ of the lift of $\gamma$ that crosses
$\til\alpha$. In Figure \ref{twist beta gamma} we indicate
$\til\gamma$, together with a lift $\til\beta$ of a representative of
$\beta$ that travels along $\til\alpha$ and the lifts of $\tau$ and
$\tau'$. The essential intersections of $\til\beta$ with $\til\gamma$
can be read off from this diagram (using the fact that all the segments
shown are geodesics) to obtain
\begin{equation}
  \label{gamma intersect tau}
 i(\til\gamma,\til\beta)-2 \le
  \#(\til\gamma\intersect(\til\tau\union\til\tau'))
\le i(\til\gamma,\til\beta)+2.
\end{equation}
(The $\pm 2$ comes from various possibilities for the endpoints and
the intersection with $\til\alpha$). Thus by 
(\ref{annulus intersection bounds distance}), 
\begin{equation}
\label{dY upper bound}
d_Y(\beta,\gamma) \le 
\#(\til\gamma\intersect(\til\tau\union\til\tau')) + 3 .
\end{equation}
Now we see from the situation in $\hhat Y$ that any two successive
intersections of $\til\gamma$ with $\til\tau$ (or $\til\tau'$) bound a segment
of $\til\gamma$ of length at least $\ell(\alpha)\ge \ep$. All these
segments project to disjoint segments in $\gamma$, and summing over
all lifts of $\gamma$ that cross $\til\alpha$, we obtain
\begin{equation}
\label{ell lower bound}
\begin{split}
\ell_{\sigma_\gamma}(\gamma) 
& \ge \ep i(\alpha,\gamma)
(\#(\til\gamma\intersect(\til\tau\union\til\tau'))-2) \\
& \ge \ep i(\alpha,\gamma)(d_Y(\beta,\gamma)-5)
\end{split}
\end{equation}
where the last line follows from (\ref{dY upper bound}).
Putting (\ref{annulus efficiency}) and (\ref{ell lower bound})
together, we find that 
\begin{equation}
d_Y(\beta,\gamma) \le 5 + (L+2C)/\ep.
\end{equation}
This is the desired bound.
\end{pf}

The analogue of Lemma \ref{pleating on common boundary} is the
following. A closely related result appears in Brock \cite{brock:continuity}.
\begin{lemma}{pleating on common curve}
For any $\ep>0$ there 
exists $D$ (depending on $\ep$ and the topological type of $S$) such
that for any Kleinian surface group
$\rho\co \pi_1(S)\to\PSL 2(\C)$ the following holds. 

Let $Y$ be an essential annulus and $g_0,g_1$ a pair of
pleated surfaces in the homotopy class $[\rho]$
mapping the core $\alpha$ of $Y$ to a 
geodesic. Let $\sigma_0$ and $\sigma_1$ be the induced metrics on $S$.
Suppose that $\ell_N(\alpha) \ge \ep>0$. 
If $\beta_0$ and $\beta_1$ are
minimal curves crossing $\alpha$ relative to  $\sigma_0$ and
$\sigma_1$, respectively, then
$$
d_Y(\beta_0,\beta_1) \le D.
$$
\end{lemma}

As in the proof of Lemma \ref{pleating on common boundary}, if
$\ell_{\sigma_0}(\tau_0) \le \delta_0$ then $\tau_0$ and $\tau_1$ are
either equal or disjoint. This bounds $i(\beta_0,\beta_1) $ and, in
turn, $d_Y(\beta_0,\beta_1)$.

If $\ell_{\sigma_0}(\tau_0)>\delta_0$ we deduce a length upper bound
on $\alpha$  as in (\ref{collar bound}), and consider the homotopy
argument of case (b) in Lemma 
\ref{pleating on common boundary}, in which we join $g_0$ to $g_1$ by
a family $g_t$. In each metric $\sigma_t$ we build $\beta_t$ as
before, and the same argument shows that only a bounded number of
curves can be built in this way. The intersection number of two such
curves that occur for the same value of $t$ is again bounded, and 
we obtain the bound on $d_Y(\beta_0,\beta_1)$ in the same way.

%
%

\section{The proof of Theorem A}
\label{proof of A}
Recall again the statement:
\state{Theorem A}{
For any Kleinian surface group $\rho$ with ending invariants
$\nu_\pm$,  if
$$
\sup_{Y\subset S} d_Y(\nu_+,\nu_-) = \infty
$$
then
$$
\inj_0(\rho) = 0,
$$
where the supremum is over proper essential subsurfaces $Y$ in $S$
not all of whose boundaries map to parabolics.
}

As in the proof of Theorem B we will work with the contrapositive
statement, that 
for any Kleinian surface group $\rho$ with ending invariants
$\nu_\pm$, 
if $\inj_0(\rho) > 0$ then
$$
\sup_{Y\subset S} d_Y(\nu_+,\nu_-) < \infty,
$$
with the supremum taken over essential subsurfaces $Y$ not all of
whose boundary components map to  parabolics.

To make sense of this statement we must first
give a general description of the ending invariants $\nu_\pm$
of a Kleinian surface group $\rho\co \pi_1(S)\to\PSL 2(\C)$.
For further
discussion see Thurston \cite{wpt:notes}, Bonahon \cite{bonahon}, and
Ohshika \cite{ohshika:ending-lams},
as well as  \cite{minsky:knoxville}.

One convenient way to think of an end invariant $\nu_+$ is 
as a pair $(\nu^G_+,\nu^T_+)$ with the following structure.
The first component $\nu^G_+$ is either 
a geodesic lamination that admits a transverse measure of full
support, or the ``empty'' lamination $\emptyset$.
The second component 
$\nu^T_+$ is either
a conformal structure of finite type on a (possibly disconnected) essential
subsurface $R_+\subset S$ (ie $\nu^T_+\in\TT(R_+)$), or
``$\emptyset$'' (in which case $R_+$ is the empty set as well).

Furthermore, $\nu^G_+$ is supported in the complement of $R_+$, and
$\nu_+$ {\em fills} in the sense that any nontrivial curve in $S$
either intersects a component of $\nu^G_+$, or of $R_+$, or it is
isotopic to a closed curve component of $\nu^G_+$.  The other end
invariant $\nu_-$ is described in the same way.

(One should keep in mind the special case $\nu^G_+=\emptyset$,
in which case $R_+=S$ and $\nu^T_+$ is the conformal structure on a
geometrically finite end of $N$, as well as the case
$\nu^T_+=\emptyset$, when
$\nu^G_+$ is a geodesic lamination that fills $S$ and the manifold has
a simply degenerate end. The other cases are hybrids of these two in
which ends of different types are separated by parabolics.)

Let $Y$ be any proper essential subsurface whose boundary components
are not all homotopic to components of $\nu^G_+$. We can define $\pi_Y(\nu_+)$
as the union of $\pi_Y(\nu^G_+)$ and $\pi_Y(\nu^T_+)$, which are
defined as follows.

$\pi_Y(\nu^G_+)$ is defined similarly to the definition of $\pi_Y$ for
closed curves: it is
the barycenter in $\CC'(Y)$ of the simplex whose vertices
are the equivalence classes of essential closed curves
and properly embedded arcs in $\nu^G_+\intersect Y$
(if there are infinite leaves of
$\nu^G_+$ that are wholly contained in $Y$, we ignore them). Note that
by choice of $Y$ this is empty only if $Y$ is contained in $R_+$.

To define $\pi_Y(\nu^T_+)$, let $\alpha_+$ be
the union of shortest curves in the conformal structure $\nu^T_+$ that
intersect $Y$ (this could be the boundary of $R_+$,
whose $\nu^T_+$-length is 0), and let $\pi_Y(\nu^T_+) =
\pi_Y(\alpha_+)$. 
This is now a
subset of $\CC'(Y)$ rather than a point, but its diameter is uniformly
bounded under one more assumption: Let $\ep>0$ be fixed and suppose
that $Y$ is anything other than an annulus whose core is homotopic to
a curve in $R_+$ with $\nu^T_+$-length  less than
$\ep$. Then (assuming $\nu^T_+\ne\emptyset$) the curves $\alpha_+$
have $\nu^T_+$--length bounded 
above by some $L(\ep)$, and hence their intersection number is bounded
above. Thus their distance in $\CC'(Y)$ is bounded above, by
Lemma \ref{bounded curves close}.

The union  $\pi_Y(\nu_+)=\pi_Y(\nu^T_+) \union \pi_Y(\nu^G_+)$ is
always nonempty under our assumptions on $Y$, 
and also has bounded diameter: Every
arc in the definition of $\pi_Y(\nu^G_+)$ is disjoint from the arcs in
$\pi_Y(\nu^T_+)$ since the supports are disjoint.

Without giving the complete definition, let us just list the
 properties of the ending invariants which we will be using.
First,   the closed
 components of $\nu^G_+$ (which include
the boundary curves of the support of $\nu^T_+$), 
 are all taken to parabolics by $\rho$. The
 same is true for $\nu_-$, and this accounts for all parabolics in
 $\rho$ other than those coming from punctures of $S$.
The Riemann surface $(R_+,\nu^T_+)$ is part of the quotient of the domain of
 discontinuity of $\rho(\pi_1(S))$, and so
by an inequality of Bers \cite{bers:boundaries},
for each curve $\gamma$ in $R_+$, 
its length $\ell_\rho(\gamma)$ is bounded above by
 $2\ell_{\nu^T_+}(\gamma)$.
Finally, for each component $\lambda$ of $\nu^G_+$ there is a sequence
 of simple closed curves $\alpha_i$ in $S$, with $\ell_\rho(\alpha_i)
 \le L_0$,  
which converge as projective measured laminations to a limit
whose support is $\lambda$.
In particular the Hausdorff limit of $\{\alpha_i\}$ contains
 $\lambda$, so that $\pi_Y(\alpha_i)$ for any fixed $Y$ is eventually
 in the same simplex as $\pi_Y(\lambda)$.

Let us now also assume that $\inj_0(\rho) = \ep>0$. Then in particular
(applying Bers' inequality from above)
the conformal structures $\nu^T_\pm$ contain no non-peripheral curves
of length less than $\ep$, and there is a bound $L_4\ge L_0$ on the
$\rho$--lengths of the curves used to construct $\pi_Y(\nu^T_\pm)$.
Putting these facts together we then have the following statement:

{\em $\pi_Y(\nu_-)$ and $\pi_Y(\nu_+)$ are contained in a
  1--neighborhood of $\pi_Y(\CC(\rho,L_4))$.}

Now the proof of the theorem is almost immediate. Let $Y$ be any
essential subsurface of $S$ not all of whose boundaries map to parabolics.
Then $\pi_Y(\nu_+)$ and $\pi_Y(\nu_-)$ are both non-empty,
and $\ell_\rho(\boundary Y) \ge \inj_0(\rho) = \ep$, and hence  by
Theorem B, $\diam_Y(\CC_0(\rho,L_4)) \le K(\ep,L_4)$.
The above observation then yields
$d_Y(\nu_-,\nu_+) \le K + 2$. This bound, independent of $Y$, gives
the desired result.

\end{document}

%% file: amsmacs.tex
\catcode`\@=11

\long\def\@savemarbox#1#2{\global\setbox#1\vtop{\hsize\marginparwidth 
  \@parboxrestore\tiny\raggedright #2}}
\newcommand\lref[1]{\ref{#1}%
\@ifundefined{r@DisplaY #1}{}{ (#1)}}

\newcommand\fakelabel[2]{\@bsphack\if@filesw {\let\thepage\relax
   \newcommand\protect{\noexpand\noexpand\noexpand}%
\xdef\@gtempa{\write\@auxout{\string
      \newlabel{#1}{{#2}{\thepage}}}}}\@gtempa
   \if@nobreak \ifvmode\nobreak\fi\fi\fi\@esphack}

\catcode`\@=12

\def\Empty{}
\newcommand\oplabel[1]{
  \def\OpArg{#1} \ifx \OpArg\Empty {} \else
        \label{#1}
  \fi}
\newtheorem{theoremSt}{Theorem}[section]

\newtheorem{exampleSt}[theoremSt]{Example}
\newtheorem{exerciseSt}[theoremSt]{Exercise}

\newcommand\MakeStEnv[1]{
  \newenvironment{#1}[1]{
  \begin{#1St} \oplabel{##1}%
  \global\def\CrntSt{\thetheoremSt}%
}{ 
  \end{#1St} }
  \newenvironment{#1+}[1]{
  \begin{#1St} \label{##1}%
  \label{DisplaY ##1}%
  \global\def\CrntSt{\thetheoremSt}%
  \def\Labl{##1}\ifx\Labl\Empty{} \else {\em (\Labl)\,}\fi%
}{ 
  \end{#1St} }
}
\MakeStEnv{theorem}
\MakeStEnv{corollary}
\MakeStEnv{proposition}
\MakeStEnv{lemma}
\MakeStEnv{definition}
\MakeStEnv{conjecture}

\long\def\state#1#2{
\medskip\par\noindent
{\bf #1 } 
{\it #2}
\par\medskip
}

\newcommand{\segment}[1]{\medskip
\begin{center} \bf #1 \end{center}
\medskip}
        
\newlength{\saveu}

\newenvironment{pf}{%
 \begin{proof}%
}{ 
 \end{proof}
}
\newenvironment{pf*}[1]{%
 \begin{proof}[#1]%
}{ 
 \end{proof}
}

\newcommand{\finishproof}[1]{ 
  \def\FPArg{#1}
  \ifx\FPArg\Empty
        \newcommand\FPArg{\CrntSt}  \fi
  \smallbreak\noindent\makebox[\textwidth]{\hfill\fbox{\FPArg}}
  \medbreak\noindent
}

\newcommand{\bfheading}[1]{\par\smallskip\noindent {\bf #1}}

\newcommand\CC{{\mathcal C}}
\newcommand\DD{{\mathcal D}}

\newcommand\FF{{\mathcal F}}

\newcommand\LL{{\mathcal L}}
\newcommand\MM{{\mathcal M}}

\newcommand\PP{{\mathcal P}}

\newcommand\TT{{\mathcal T}}

\newcommand\PMF{{\PP\kern-2pt\MM\FF}}

\newcommand\PML{{\PP\kern-2pt\MM\LL}}

\newcommand\ep{\epsilon}

\newcommand\hhat{\widehat}

\newcommand\union{\cup}
\newcommand\intersect{\cap}
\newcommand\bbR{{\mathord{\text{I\kern-2pt R}}}}        
\newcommand\bbH{{\mathord{\text{I\kern-2pt H}}}}        

\newcommand\C{{\mathbf C}}

\newcommand\Z{{\mathbf Z}}
\newcommand\R{{\mathbf R}}
\newcommand\Q{{\mathbf Q}}

\newcommand\Hyp{{\mathbf H}}

\newcommand\PSL[1]{\text{PSL}_{#1}}

\newcommand\bigrightarrow[1]{\hbox to #1{\rightarrowfill}}
\newcommand\bigleftarrow[1]{\hbox to #1{\leftarrowfill}}

\newcommand\boundary{\partial}
\newcommand\semidir{\mathrel{\hbox{\vrule depth-.03ex height1.1ex\kern-0.15em$\times$}}}

\newcommand\til{\widetilde}

\newcommand{\diam}{\operatorname{diam}}

\numberwithin{equation}{section}